\documentclass{article}
\usepackage[utf8]{inputenc}
\usepackage{graphicx}
\usepackage{float}
\usepackage{amsmath}
\usepackage{footnotehyper}
\usepackage[colorlinks,hyperindex=true,bookmarks=true]{hyperref}

\newtheorem{theorem}{Theorem}

\author{Kobelkov~S.G.\footnote{Lomonosov Moscow State University,Moscow,Russia}\footnote{HSE,Moscow,Russia}, Mishchenko~A.A.\footnote{LLC TD Uralprotect}}

\def\calR{\mathcal{R}}
\def\grad{\nabla}

\begin{document}
\title{Irreversible Step-Growth Polymerization with weighted dynamic networks}
\abstract{
	This paper considers random graph approach to simulate irreversible step-growth polymerization. We study generalization of approach developed by Kryven~I. by introducing different types of bonds each with its own weight representing corresponding reaction rate. 
	%A special case of dynamic functional degree change is also considered.
}
\maketitle

	\section{Introduction}
	
	Many research fields benefit from recent developments in random graphs and network theory. In particular, polymerization process can be described well with random graph models \cite{Newman,Kryven42,Kryven} and these results are consistent with well-known models, e.g. Flory \cite{Flory} or Stockmayer \cite{Stockmayer}. 
	
	Conventional polymer networks are formed by a process called polymerization during which monomers bind together by means of covalent bonds and form large molecular structures. One of the most common polymerization processes is the step-growth polymerization of multifunctional monomers. This process leads to hyperbranched polymers of irregular topologies that undergo a phase transition in their connectivity structure which is marked by emergence of the gel, i.e. the giant molecule that spans the whole volume. 
	Flory provided simple analytical expressions for the average molecular size for several particular configurations. Later Stockmayer presented a formal expression for the whole distribution of molecular sizes, but its application is quite limited because of complicated computations. 
	
	Here we should mention fast approximated methods as in papers \cite{Wulkow,Tobita,Hillegers,Kryven35}. However, these are very specific and hard to adapt to new polymerization schemes. Molecular dynamic simulation approaches provide lots of details on the structure of polymer networks, however they are computationally expensive and, thus, limited to small samples and short time scales.
	
    On the other hand random graph approaches present a opportunity to obtain exactly solvable expressions for hyperbranched polymers.
	In the case of symmetric covalent bonds the structure of polymer networks is described by the configuration model for undirected random networks \cite{Newman,Kryven42}. Developments in directed configuration models allowed to develop a generic polymerization framework that covers asymmetrical bonds as well. 
	
	Unfortunately, these models miss the case when both symmetrical and asymmetrical bonds are present in the system, and do not account for rate at which these bonds are formed. This paper tries to address these issues.

	\section{Master equation}
	
	We start from considering simple case when we have two types of monomers, say, A and B. A can connect with B, but not with self; on the opposite, B can connect both with A and B. This differs from the models considered in other papers by mixing both directional and non-directional connections.

We distinguish monomer species by counting the numbers of functional groups
of both types $I$, $J$ and the numbers of in- and bidirectional-edges $i, j$. During the progress of polymerisation the functional
groups are converted into chemical bonds between the monomers, and the concentration profiles $M_{i,j,I,J}(t)$ obey
the following master equation:	
	\begin{equation}\label{mastereq}
		\begin{array}{l}
		\frac{\partial}{\partial t}M_{i,j,I,J}(t)=(I-i+1)(\nu_{01}-\mu_{01}(t))M_{i-1,j,I,J}(t)+\\
		(J-j+1)(\nu_{10}-\mu_{10}(t)+\nu_{01}-\mu_{01}(t))M_{i,j-1,I,J}(t)
		-\\
		((I-i)(\nu_{01}-\mu_{01}(t))+(J-j)(\nu_{10}-\mu_{10}(t)+\nu_{01}-\mu_{01}(t)))M_{i,j,I,J}(t)
				\end{array}		
	\end{equation}
Like in the case of fully oriented graph, considered in \cite{Kryven}, probability of connecting functional group A is proportional to the number of free functional groups B. On the other hand, since the functional groups B can connect both to A and B, probability of connecting them is proportional to the total number of free functional groups.

Initially at $t=0$ there are no bonds, thus we have the following 
initial condition 
\begin{equation}\label{initcond}
\begin{array}{ll}
M_{i,j,I,J}(0)=P(I,J) & \mbox{if } i=j=0,\\
M_{i,j,I,J}(0)=0 & \mbox{otherwise.}
\end{array}
\end{equation}
Here $P(I,J)$ denotes initial distribution of monomers of $I,J$ type. 
Define $\nu_{mn}$ to be partial moments of the initial monomer type distribution
\[
\nu_{mn}=\sum_{I,J\ge 0}I^mJ^nP(I,J)
\]

Let the time dependent degree distribution $u(i, j, t)=\sum_{I,J}M_{i,j,I,J}(t)$, then 
\begin{equation}\label{mu_mn}
\mu_{mn}(t)=\sum_{i,j\ge 0}i^mj^nu(i,j,t),
\end{equation}

e.g. $\mu_{mn}$ are partial moments of degree distribution. Note that above definition is not strict since summation is performed up to $I$ and $J$ correspondingly, so one should substitute expression for $u(i,j,t)$ and change summation order, e.g. $\mu_{mn}(t)=\sum_{I,J}\sum_{i,j}i^mj^nM_{i,j,I,J}(t)$.

In order to solve the system of differential equations (\ref{mastereq},\ref{initcond}) we proceed to generating functions, namely,
let $m(x,y,I,J,t)=\sum_{i=0}^I\sum_{j=0}^J x^iy^jM_{i,j,I,J}(t)$. 

Obviously, we have the following relations:
\begin{equation}\label{relations_m}
\begin{array}{l}
	\frac{\partial}{\partial x}m(x,y,I,J,t)=\sum_{i\ge 1,j}ix^{i-1}y^jM_{i,j,I,J}(t),\\
	\frac{\partial}{\partial y}m(x,y,I,J,t)=\sum_{j\ge 1,i}jx^{i}y^{j-1}M_{i,j,I,J}(t),\\
	m(x,y,I,J,0)=P(I,J).
\end{array}
\end{equation}

Let us multiply equations (\ref{mastereq}) by $x^iy^j$ and sum up over all $i,j$. Taking advantage of (\ref{relations_m}), we obtain the following system of PDEs
\begin{equation}\label{pde_m}
	\begin{array}{l}
\frac{\partial m(x,y,I,J,t)}{\partial t}=(\nu_{01}-\mu_{01}(t))(Ixm(x,y,I,J,t)-x^2\frac{\partial m(x,y,I,J,t)}{\partial x})+\\
(\nu_{10}-\mu_{10}(t)+\nu_{01}-\mu_{01}(t))(Jym(x,y,I,J,t)-y^2\frac{\partial m(x,y,I,J,t)}{\partial y})-\\
(\nu_{01}-\mu_{01}(t))(Im(x,y,I,J,t)-x\frac{\partial m(x,y,I,J,t)}{\partial x})-\\
(\nu_{01}-\mu_{01}(t)+\mu_{10}-\mu_{10}(t))(J(x,y,I,J,t)m-y\frac{\partial m(x,y,I,J,t)}{y}),\\
m(x,y,I,J,0)=P(I,J).
	\end{array}
\end{equation}

Since $u(i,j,t)=\sum_{I,J}M_{i,j,I,J}(t)$, we can write
\begin{equation}\label{mu_through_m}
	\begin{array}{l}
\mu_{10}(t)=\sum_{I,J}\left.\frac{\partial m(x,y,I,J,t)}{\partial x}\right|_{x=1,y=1},\\
\mu_{01}(t)=\sum_{I,J}\left.\frac{\partial m(x,y,I,J,t)}{\partial y}\right|_{x=1,y=1}.
	\end{array}
\end{equation}

Next, we're going to obtain a system of ODEs for the functions $\mu_{10}(t)$ and $\mu_{01}(t)$.

Let us rewrite the differential equation from (\ref{pde_m}) as
\begin{equation}\label{m_short}
	\begin{array}{l}
\frac{\partial m(x,y,I,J,t)}{\partial t}=(\nu_{01}-\mu_{01}(t))(x-1)(Im(x,y,I,J,t)-x\frac{\partial m(x,y,I,J,t)}{\partial x})+\\
(\nu_{01}-\mu_{10}(t)+\mu_{01}-\mu_{01}(t))(y-1)(Jm(x,y,I,J,t)-y\frac{\partial m(x,y,I,J,t)}{\partial y})
\end{array}
\end{equation}
Differentiating (\ref{m_short}) by $x$ and taking $x=1,y=1$, we have
\[
\left.\frac{\partial m(x,y,I,J,t)}{\partial t\partial x}\right|_{x=y=1}=(\nu_{01}-\mu_{01}(t))\left(Im(1,1,I,J,t)-\left.\frac{\partial m(x,y,I,J,t)}{\partial x}\right|_{x=y=1}\right)
\]
After summing up over all $I,J$ and taking into account (\ref{mu_through_m}) we have
\[
\frac{\partial}{\partial t}\mu_{10}(t)=(\nu_{01}-\mu_{01}(t))(\sum_{I,J}Im(1,1,I,J,t)-\mu_{10}(t)).
\]
Note that $\sum_{I,J}Im(1,1,I,J,t)=\sum_{I,J}I\sum_{i,j}M_{i,j,I,J}(t)$. Since total number of monomer species of each type remains the same, $\sum_{i,j}M_{i,j,I,J}(t)=\sum_{i,j}M_{i,j,I,J}(0)$, thus
$\sum_{I,J}Im(1,1,I,J,t)=\sum_{I,J}IP(I,J)=\nu_{10}$. 

Finally, $\mu_{10}(0)=\sum_{I,J}\left.\frac{\partial m(x,y,I,J,t)}{\partial x}\right|_{x=y=1,t=0}$ and it remains to recall that $m(x,y,I,J,0)=P(I,J)$ for all $x,y$, so $\mu_{10}(0)=0$.

Differentiating (\ref{m_short}) by $y$, taking $x=1,y=1$, summing up by $I,J$ and taking into account (\ref{mu_through_m}), we obtain system of ODEs for $\mu_{01}(t)$ and $\mu_{10}(t)$.

\begin{equation}\label{ode_mu}
	\begin{array}{l}
	\frac{\partial\mu_{10}(t)}{\partial t}=(\nu_{01}-\mu_{01}(t))(\nu_{10}-\mu_{10}(t)),\\
	\frac{\partial\mu_{01}(t)}{\partial t}=(\nu_{01}-\mu_{01}(t)+\nu_{10}-\mu_{10}(t))(\nu_{01}-\mu_{01}(t)),\\
	\mu_{01}(0)=\mu_{10}(0)=0.
	\end{array}
\end{equation}

In \cite{Kryven2} they consider fully oriented graph, which makes the system (\ref{ode_mu}) symmetrical, e.g. the second equation looks like
$\frac{\partial\mu_{01}(t)}{t}=(\nu_{10}-\mu_{10}(t))(\nu_{01}-\mu_{01}(t))$; evidently, one can conclude $\mu_{01}(t)=\mu_{10}(t)$ and can easily find analytical solution of such a system. 

In our case $\mu_{01}(t)\not\equiv\mu_{10}(t)$, so we're solving (\ref{ode_mu}) numerically with RK methods.
%For example, if we take $\nu_{01}=\nu_{10}=1$, solution of the system can be expressed through solution of the equation
%\[
%\int_0^x \frac{1}{(z-1)^2(1-\ln(z))}dz=t
%\]
%which is known not to have closed form.

%Moreover, numerical solutions for different values of $\nu_{01}$ and $\nu_{10}$ show that there exists $t_0>0$, such that $\mu_{10}(t_0)=0$ and $\mu_{10}(t)<0$ for $t>t_0$; finally, there is a singularity point $t>t_0$ where solution is not defined. We believe that $t_0$ corresponds to the moment where there are no free functional groups of type B left, thus polymerization process stops.

Now, let's go back to solving the system (\ref{pde_m}). Assume solution in form 
\[
m(x,y,I,J,t)=\left(\frac{1+A(t)t}{1+A(t)tx}\right)^{-I}\left(\frac{1+B(t)t}{1+B(t)ty}\right)^{-J}P(I,J).
\]
Consequently,
\begin{equation}\label{log_m}
	\begin{array}{l}
	\log m(x,y,I,J,t)=-I\log(1+A(t)t)+I\log(1+A(t)tx)-\\
	J\log(1+B(t)t)+J\log(1+B(t)ty)+\log P(I,J).
	\end{array}
\end{equation}

Let us rewrite (\ref{m_short}), dividing left and right parts by $m(x,y,I,J,t)$:
\[
\begin{array}{l}
\frac{\partial\log m(x,y,I,J,t)}{\partial t}=(\nu_{01}-\mu_{01}(t))(x-1)(I-x\frac{\partial\log m(x,y,I,J,t)}{\partial x})+\\
(\nu_{10}-\mu_{10}(t)+\nu_{01}-\mu_{01}(t))(y-1)(J-y\frac{\partial\log m(x,y,I,J,t)}{\partial y})
\end{array}
\]
and substitute (\ref{log_m}). 

Equate the terms depending on $A(t)$:
\[
\frac{I(A(t)+t\frac{\partial A(t)}{\partial t})(x-1)}{(1+A(t)t)(1+A(t)tx)}=
\frac{I(\nu_{01}-\mu_{01}(t))(x-1)}{1+A(t)tx}
\]
Thus,
\[
\frac{A(t)+t\frac{\partial A(t)}{\partial t}}{1+A(t)t}=\nu_{01}-\mu_{01}(t)
\]
with initial condition is $A(0)=\nu_{01}$. Consequently,
\begin{equation}\label{At}
A(t)=\frac{\exp\{\int_0^t \nu_{01}-\mu_{01}(u)du\}-1}{t}.
\end{equation}

Similarly,
\begin{equation}\label{Bt}
	B(t)=\frac{\exp\{\int_0^t \nu_{10}-\mu_{10}(u)+\nu_{01}-\mu_{01}(u)du\}-1}{t}.
\end{equation}

Since $\sum_{i=0}^n\binom{n}{k}x^ip^i(1+p)^{-n}=\left(\frac{1+p}{1+xp}\right)^{-n}$, we get expression for $M_{i,j,I,J}(t)$:
\[
M_{i,j,I,J}(t)=\binom{I}{i}(tA(t))^i(1+tA(t))^{-I}\binom{J}{j}(tB(t))^j(1+tB(t))^{-J}P(I,J).
\]
This expression can be rewritten using binominal probabilities. Let us denote $p_A(t)=\frac{tA(t)}{1+tA(t)}$ and
$p_B(t)=\frac{tB(t)}{1+tB(t)}$, then
\begin{equation}\label{expr_M}
M_{i,j,I,J}(t)=\binom{I}{i}p_A(t)^i(1-p_A(t))^{I-i}\binom{J}{j}p_B(t)^j(1-p_B(t))^{J-j}P(I,J).
\end{equation}

Consequently,
\begin{equation}\label{expr_u}
	u(i,j,t)=\sum_{I,J}\binom{I}{i}p_A(t)^i(1-p_A(t))^{I-i}\binom{J}{j}p_B(t)^j(1-p_B(t))^{J-j}P(I,J)
\end{equation}

Once expression for $u(i,j,t)$ is set up, we can calculate $\mu_{nm}(t)$ according to (\ref{mu_mn}), namely
\begin{equation}\label{mu_mn_e}
	\begin{array}{l}
\mu_{00}(t)=1,\\
\mu_{10}(t)=p_A(t)\nu_{10},\mu_{01}(t)=p_B(t)\nu_{01},\\
\mu_{20}(t)=p_A(t)(p_A(t)\nu_{20}-p_A(t)\nu_{10}+\nu_{10}),\\
\mu_{02}(t)=p_B(t)(p_B(t)\nu_{02}-p_B(t)\nu_{01}+\nu_{01}),\\
\mu_{11}(t)=p_A(t)p_B(t)\nu_{11},
\mu_{22}(t)=p_A(t)p_B(t)(p_A(t)p_B(t)(\nu_{22}-\nu_{21}-\nu_{12}+\nu_{11})+\\
p_A(t)(\nu_{21}-\nu_{11})+p_B(t)(\nu_{12}-\nu_{11})+\nu_{11}).
	\end{array}
\end{equation}

Note that obtained relations for $\mu_{10}(t)$ and $\mu_{01}(t)$ allow much faster calculation of $p_A(t)$ and $p_B(t)$ compared to their definitions.
Stationary solution for the system (\ref{ode_mu}) can be partially obtained by equating left-hand sides to zero, thus we immediately have $\mu_{01}(t)\rightarrow\nu_{01}$ as $t\rightarrow\infty$.

\textbf{Systems with a single monomer type}

Consider the case when there is only one monomer species bearing $I$ groups of type $A$ and $J$ groups of type $B$.
This corresponds to condition $P(I,J)=1$, thus all the sums over all $I$ and $J$ reduce to a single summand.
Moments of distribution $P$ are $\nu_{mn}=I^mJ^n$. 

Unfortunatelly, even in this case the system (\ref{ode_mu}) cannot be solved analytically, and we should stick to numerical solution.
%\begin{figure}
%	\centering
%%	\includegraphics[width=0.5\linewidth]{sol11}
%	\caption{$\mu_{01}$ and $\mu_{10}$ for $I=1,J=1$}
%	\label{fig:sol11}
%\end{figure}
%
%\begin{figure}
%	\centering
%%	\includegraphics[width=0.5\linewidth]{sol44}
%	\caption{$\mu_{01}$ and $\mu_{10}$ for $I=4,J=4$}
%	\label{fig:sol44}
%\end{figure}
%
%\begin{figure}
%	\centering
%%	\includegraphics[width=0.5\linewidth]{sol2g4f}
%	\caption{$\mu_{01}$ and $\mu_{10}$ for $I=4,J=2$}
%	\label{fig:sol2g4f}
%\end{figure}
%
%\begin{figure}
%	\centering
%%	\includegraphics[width=0.5\linewidth]{sol4g2f}
%	\caption{$\mu_{01}$ and $\mu_{10}$ for $I=2,J=4$}
%	\label{fig:sol4g2f}
%\end{figure}
%

Consequently we can calculate $A(t),B(t)$, and $p_A(t),p_B(t)$. Numerical results for this kind of model are presented in Appendix.

%
%\begin{figure}
%	\centering
%%	\includegraphics[width=0.5\linewidth]{pApB11}
%	\caption{$p_A(t)$ and $p_B(t)$ for $I=1,J=1$}
%	\label{fig:pApB1_1}
%\end{figure}
%
%\begin{figure}
%	\centering
%%	\includegraphics[width=0.5\linewidth]{pApB44}
%	\caption{$p_A(t)$ and $p_B(t)$ for $I=4,J=4$}
%	\label{fig:pApB4_4}
%\end{figure}
%
%\begin{figure}
%	\centering
%%	\includegraphics[width=0.5\linewidth]{pApB4_2}
%	\caption{$p_A(t)$ and $p_B(t)$ for $I=4,J=2$}
%	\label{fig:pApB4_2}
%\end{figure}
%
%\begin{figure}
%	\centering
%%	\includegraphics[width=0.5\linewidth]{pApB2_4}
%	\caption{$p_A(t)$ and $p_B(t)$ for $I=2,J=4$}
%	\label{fig:pApB2_4}
%\end{figure}

\section{Gel transition point}

Next, we'll try to derive gen transition point.
Let $U(x,y,t)$ be generating function for degree distribution $u(i,j,t)$. Suppose, we select a vertex that is at the end of a directed randomly-chosen edge. Conditional degree distribution in this case is 
$u_{in}(n,k,t)=\frac{n}{\mu_{10}(t)}u(n,k,t)$. Correspondingly, if we select a vertex at non-directed randomly-chosen edge, conditional degree distribution would be $u_{bi}(n,k,t)=\frac{k}{\mu_{01}(t)}u(n,k,t)$. Corresponding generating functions
can be expressed as 
\begin{equation}\label{def_U_in_U_bi}
\begin{array}{l}
U_{in}(x,y,t)=\frac{1}{\mu_{10}(t)}\frac{\partial U(x,y,t)}{\partial x},\\
U_{bi}(x,y,t)=\frac{1}{\mu_{01}(t)}\frac{\partial U(x,y,t)}{\partial y}.
\end{array}
\end{equation}

If $w(s,t)$ is the probability that a randomly chosen monomer belongs to a component of size $s$ (e.g. molecular weight distribution), then
$w(1,t)=u(0,0,t), w(2,t)=u(1,0,t)u(0,1,t)\left(\frac{1}{\mu_{01}(t)}+\frac{1}{\mu_{10}(t)}\right)+\frac{u(0,1,t)^2}{\mu_{01}(t)}$, etc.

Let $W(x,t)$ be GF of $w(s,t)$. Now consider biased choice for the starting vertex. Suppose, one selects a directed edge at random and the end vertex of this edge as root. Then, $w_{in}(s,t)$ denotes the distribution of weak-component sizes associated with this root; corresponding GF we denote through $W_{in}(x,t)$. Similarly, for the non-directed edge we define GF $W_{bi}(x,t)$. Bind $W_{in}, W_{bi}, W, U_{in}, U_{bi}, U$ together.

Let us start by selecting a root vertex that we arrive at by following a directed edge. The probability of the root to have $n$ in-edges and $k$ non-directed edges is $u_{in}(n,k,t)$. Each of the in-edges is associated with a weak component of the size $w_{bi}(s,t)$, thus the sum of sizes of all components reached through the in-edges is distributed as a sum of independent copies of $w_{bi}(s,t)$. On the other hand, each non-directional edge can be associated with a weak component of the size $w_{in}(s,t)$ with conditional probability proportional to $\mu_{10}(t)$ and with a weak component of the size $w_{bi}(s,t)$ with probability proportional to $\mu_{01}(t)$. Let $\kappa(t)=\frac{\mu_{10}(t)}{\mu_{10}(t)+\mu_{01}(t)}$.
Introduce $W_{+}(s,t)=\kappa(t)W_{in}(x,t)+(1-\kappa(t))W_{bi}(x,t)$, then the generating function for a weak component distrubution reached through the non-directional edge is $W_{+}(s,t)$.

Thus, the distribution for the sum of sizes of all components originated at the root is obtained as 
\[
\sum_{n,k}u_{in}(n,k,t)W_{bi}(x,t)^nW_{+}(x,t)^k=U_{in}(W_{bi}(x,t),W_{+}(x,t)).
\]
On the other hand, the total number of all vertices reachable from the root plus one can be also considered as the size of the weak component reached by following the original in-edge. Thus, one obtains the following relation:
\[
W_{in}(x)=xU_{in}(W_{bi}(x,t),W_{+}(x,t)).
\]
Similarly, we get the following system of functional equations:
\begin{equation}\nonumber
	\begin{array}{l}
		W(x,t)=xU(W_{bi}(x,t),W_{+}(x,t),t),\\
		W_{in}(x,t)=xU_{in}(W_{bi}(x,t),W_{+}(x,t),t),\\
		W_{bi}(x,t)=xU_{bi}(W_{bi}(x,t),W_{+}(x,t),t),
	\end{array}
\end{equation}
or equivalently,
\begin{equation}\label{W_in_bi}
	\begin{array}{l}
		W(x,t)=xU(W_{bi}(x,t),W_{+}(x,t),t),\\
		W_{+}(x,t)=xU_{+}(W_{bi}(x,t),W_{+}(x,t),t),\\
		W_{bi}(x,t)=xU_{bi}(W_{bi}(x,t),W_{+}(x,t),t),
\end{array}
\end{equation}
where $U_{+}(n,k,t)=\kappa(t)U_{in}(n,k,t)+(1-\kappa)U_{bi}(n,k,t)$.

Weight-average molecular weight is calculated as 
$\frac{W'(1,t)}{W(1)}$.

Since $\mu_{mn}(t)=\left.\left(x\frac{\partial}{\partial x}\right)^m\left(y\frac{\partial}{\partial y}\right)^nU(x,y,t)\right|_{x=y=1}$, we have $\mu_{10}(t)=U_{in}(1,1,t)\mu_{10}(t)$, e.g. $U_{in}(1,1,t)=1$. Similarly we have $U_{bi}(1,1,t)=1$. 

Consequently, (\ref{W_in_bi}) gives us $W_{in}(1,t)=1$ and $W_{bi}(1,t)=1$. It remains to calculate derivative of $W(x,t)$. Let us take derivative in $x$ for all relations in (\ref{W_in_bi}) at $x=1$:
\[
\begin{array}{l}
	\frac{\partial W(1,t)}{\partial x}=1+\frac{\partial U(1,1,t)}{\partial x}\frac{\partial W_{bi}(1,t)}{\partial x}+\frac{\partial U(1,1,t)}{\partial y}\frac{\partial W_{+}(1,t)}{\partial x}\\
	\frac{\partial W_{+}(1,t)}{\partial x}=1+\frac{\partial U_{+}(1,1,t)}{\partial x}\frac{\partial W_{bi}(1,t)}{\partial x}+\frac{\partial U_{+}(1,1,t)}{\partial y}\frac{\partial W_{+}(1,t)}{\partial x}\\
	\frac{\partial W_{bi}(1,t)}{\partial x}=1+\frac{\partial U_{bi}(1,1,t)}{\partial x}\frac{\partial W_{bi}(1,t)}{\partial x}+\frac{\partial U_{bi}(1,1,t)}{\partial y}\frac{\partial W_{+}(1,t)}{\partial x}.
\end{array}
\]

Solution to this linear system degenerates if 
\begin{equation}\label{sys_deg_cond}
	\frac{\partial U_{+}(1,1,t)}{\partial x}\frac{\partial U_{bi}(1,1,t)}{\partial y}=(\frac{\partial U_{+}(1,1,t)}{\partial y}-1)(\frac{\partial U_{bi}(1,1,t)}{\partial x}-1)
\end{equation}
and corresponds to unlimited growth of the weak component.

From (\ref{def_U_in_U_bi}) we obtain
\begin{equation}\label{diff_u_through_m}
\begin{array}{l}
	\frac{\partial U_{bi}(1,1,t)}{\partial y}=\frac{\mu_{02}(t)}{\mu_{01}(t)}-1,\\
	\frac{\partial U_{bi}(1,1,t)}{\partial x}=\frac{\mu_{11}(t)}{\mu_{01}(t)},\\
	\frac{\partial U_{in}(1,1,t)}{\partial x}=\frac{\mu_{20}(t)}{\mu_{10}(t)}-1,\\
	\frac{\partial U_{in}(1,1,t)}{\partial y}=\frac{\mu_{11}(t)}{\mu_{10}(t)}.
\end{array}
\end{equation}
Substituting this in its turn into the relation (\ref{sys_deg_cond}), we get the criteria for gel transition point.

\begin{theorem}
	Systems with two types of functional groups have gel transition point defined by the following equation
\begin{equation}\label{gel_transition}
    \begin{array}{l}
    	((\nu_{01}-\nu_{02})(\kappa(t)(1+p_A(t))-1)\nu_{10}-
    	((\nu_{01}\nu_{20}-\nu_{02}\nu_{20}+\\
    	\nu_{11}^2)p_A(t)-\nu_{01}\nu_{11})\kappa(t))p_B(t)+\nu_{10}(\nu_{11}p_A(t)-\nu_{01})=0
    \end{array}
\end{equation}
\end{theorem}

We'll use this statement to obtain numerical results in the Appendix.

%Solving the obtained system, we get $\frac{\partial W(1,t)}{\partial x}=1+\frac{A_1(t)}{B_1(t)}$, where $A_1(t)$ is continuous function, and $B_1(t)=(\nu_{11}^2+(\nu_{01}-\nu_{02})p_A(t)^2p_B(t)^2-p_B(t)p_A(t)\nu_{11}(\mu_{01}(t)+\mu_{10}(t))+\mu_{01}(t)\mu_{10}(t)$.
%
%Note that if $B_1(t)=0$, weight-average molecular weight is exploding, thus we consider this as gel transition point.
%
%Solving $B_1(t)=0$ for $p_A(t)p_B(t)$ gives the following criteria for the gel transition:
%\begin{equation}\label{p_crit}
%	\begin{array}{l}
%  p_A(t)p_B(t)=\\
%  \frac{-\nu_{11}(\mu_{01}(t)+\mu_{10}(t))+\sqrt{\nu_{11}^2(\mu_{01}(t)-\mu_{10}(t))+
%4\mu_{01}(t)\mu_{10}(t)(\nu_{10}-\nu_{20})(\nu_{01}-\nu_{02})}}{
%-2\nu_{11}^2+2(\nu_{01}-\nu_{02})(\nu_{10}-\nu_{20})
%}
%\end{array}
%\end{equation}
%
%Note that this relation is time-dependent, however, if we take a symmetrical case studied in Kryven, e.g. $\mu_{01}(t)=\mu_{10}(t)$, we can reduce this to a quite simple condition, not depending on $t$ in the right-hand side (see Kryven).

\subsection{Systems with two monomer types}

Now, let us proceeed to the other commonly used case. We consider one monomer type to handle only groups of type A, while the other contains groups of type B only. In terms of initial distribution, we can express this as 
$P(I,0)=\alpha, P(0,J)=1-\alpha$, where $\alpha$ is fraction of monomers of the first kind. Consequently, $\nu_{10}=I\alpha,\nu_{01}=J(1-\alpha),\nu_{11}=0$,$\nu_{20}=I^2\alpha,\nu_{02}=J^2(1-\alpha)$.

Like in the previous section, introduce $U(x,y,t)$ be generating function for $u(i,j,t)=\binom{I}{i}p_A(t)^i(1-p_A(t))^{(I-i)}P(I,0)\delta_j+\binom{J}{j}p_B(t)^j(1-p_B(t))^{(J-j)}P(0,J)\delta_i$.
Here $\delta_i$ is a $\delta$-function, e.g. it equals zero unless  $i=0$.

Substituting into (\ref{def_U_in_U_bi}) gives
\begin{equation*}
	\begin{array}{l}
		U(x,y,t)=((x-1)p_A(t)+1)^IP(I,0)+((y-1)p_B(t)+1)^JP(0,J),\\
		U_{in}(x,y,t)=\frac{P(I,0)}{\mu_{10}(t)}Ip_A(t)((x-1)p_A(t)+1)^{I-1},\\
		U_{bi}(x,y,t)=\frac{P(0,J)}{\mu_{01}(t)}Jp_B(t)((y-1)p_B(t)+1)^{J-1}.
	\end{array}
\end{equation*}

Since $\nu_{11}=0$ we can simplify expression (\ref{gel_transition}) for the gel transition point
\[
    	((\nu_{01}-\nu_{02})(\kappa(t)(1+p_A(t))-1)\nu_{10}-
(\nu_{01}\nu_{20}-\nu_{02}\nu_{20})p_A(t)\kappa(t))p_B(t)-\nu_{10}\nu_{01}=0
\]

Numerical results corresponding to this case are presented in the Appendix.

\subsection{Systems with three monomer types}

Considerations in this section are the same as in previous ones. The only change is the initial setup, where fractions need to be distributed between three types, say $P(2,0)=\alpha,P(0,2)=\beta,P(0,5)=1-\alpha-\beta$, $\alpha>0,\beta>0,\alpha+\beta<1$. In this case $\nu_{10}=2\alpha,\nu_{20}=4\alpha,\nu_{01}=2\beta+5(1-\alpha-\beta),\nu_{02}=4\beta+25(1-\alpha-\beta),\nu_{11}=0$.

Numerical results corresponding to this case are presented in the Appendix.

%\begin{figure}
%	\centering
%%	\includegraphics[width=0.7\linewidth]{surf_crit_2_2_5}
%	\caption{$p_{A,crit}$ surface on $\alpha,\beta$}
%	\label{fig:surfcrit225}
%\end{figure}
%

\section{Weighted polymerization}

Let us generalize the considered model by introducing weights which express reaction rate. Namely, let us have two types of monomers, $A$ and $B$. As before, $A$ can connect with $B$, but not with self; however, $B$ can react  not only with $A$, but also with $B$, and let us account that the latter is less probable. We will express this change by introducing the weight $0\le w\le 1$ which denotes ``likelyness'' of establishing connection $B-B$. 

First of all, let us see how the original master equation (\ref{mastereq}) changes:
\begin{equation}\label{master_eq_weighted}
	\begin{array}{l}
		\frac{\partial}{\partial t}M_{i,j,I,J}(t)=(I-i+1)(\nu_{01}-\mu_{01}(t))M_{i-1,j,I,J}(t)+\\
		(J-j+1)(\nu_{10}-\mu_{10}(t)+w(\nu_{01}-\mu_{01}(t)))M_{i,j-1,I,J}(t)
		-\\
		((I-i)(\nu_{01}-\mu_{01}(t))+(J-j)(\nu_{10}-\mu_{10}(t)+w(\nu_{01}-\mu_{01}(t))))M_{i,j,I,J}(t)
		\end{array}
\end{equation}
while initial conditions (\ref{initcond}) remain the same. Consequently, we proceed to GFs and obtain the following system of ODEs which resembles (\ref{ode_mu}):
\begin{equation}\label{ode_mu_weighted}
	\begin{array}{l}
		\frac{\partial\mu_{10}(t)}{t}=(\nu_{01}-\mu_{01}(t))(\nu_{10}-\mu_{10}(t)),\\
		\frac{\partial\mu_{01}(t)}{t}=(w(\nu_{01}-\mu_{01}(t))+\nu_{10}-\mu_{10}(t))(\nu_{01}-\mu_{01}(t)),\\
		\mu_{01}(0)=\mu_{10}(0)=0.
	\end{array}
\end{equation}

Once we have a numerical solution for (\ref{ode_mu_weighted}), we can assume (\ref{log_m}), where
$A(t)$ is calculated as before (\ref{At}), but 
\begin{equation*}\label{Bt_weighted}
	B(t)=\frac{\exp\{\int_0^t \nu_{10}-\mu_{10}(u)+w(\nu_{01}-\mu_{01}(u)du)\}-1}{t}.
\end{equation*}

Further considerations remain alsmost the same. If we take $\kappa(t)=\frac{\mu_{10}(t)}{\mu_{10}(t)+w\mu_{01}(t)}$ we can proceed to evaluating gel transition points depending on the weight using the same relation (\ref{gel_transition}).

Numerical results evaluating dependency of the gen transition point on weight for the case of two monomers are presented in the Appendix.

%
%\begin{figure}
%	\centering
%%	\includegraphics[width=0.7\linewidth]{plot_crit_I2J5w07}
%	\caption{Criteria (\ref{weak_giant_component}) for $w=0.7; \alpha=0.5$. $p_{A,crit}=0.2601; p_{B,crit}=0.9613$}
%	\label{fig:plotcriti2j5w07}
%\end{figure}
%
%Need to note that lower values of $w$ make it slower until gel transition point, so one should increase the simulation time frame

\section{Systems with multiple functional group types}

In this section we consider generalization of the approach to  monomer species which can have functional groups of types $I_1,I_2,\ldots,I_r, r\ge 2$. Let us define a non-degenerate non-negative symmetric weight matrix $W\in\calR^{r\times r}$, so that 
each element $w_{mn}$ defines relative rate at which $I_m,I_n$ can connect. Assume that sum of values for each row of the matrix equals $1$; zero values mean that corresponding connection is not possible.

Let further $\vec{I}$ denote vector $(I_1,I_2,\ldots,I_r)$, $\vec{\nu}_1$ be the first moment of initial distribution, e.g. ${\nu}_{1,k}=\sum_{I_k}I_kP(\vec I)$, where $P(\vec I)$ stands for initial distribution of the monomer species of different types. Also define second moment matrix $\nu_2$ with $\nu_{2,mn}=\sum_{I_m,I_n}I_mI_nP(\vec I)$.

As before, we start from writing the master equation for the $k$-th component of $M_{\vec i,\vec I}(t)$:
\begin{equation}\label{master_eq3}
	\begin{array}{l}
	\frac{\partial}{\partial t}M_{\vec i,\vec I}(t)=
	((\vec I-\vec i+1)\odot M_{\vec i-\vec e,\vec I(t)}(t)-(\vec I-\vec i)M_{\vec i,\vec I}(t))^TW(\vec\nu_1-\vec\mu),
	\end{array}
\end{equation}
where $M_{\vec i-\vec e,\vec I}(t)$ denotes a vector with components $M_{\vec i-\vec e_1,\vec I}(t)$, and $\vec x\odot\vec y$ stands for Hadamard product.
Initial condition is $M_{\vec i,\vec I}(0)=P(\vec I)$ if $\vec i=0$ and zero otherwise. 

Let the time dependent degree distribution $u(\vec i,t)=\sum_{\vec I}M_{\vec i,\vec I}(t)$, then
\[
\mu_{m}(t)=\sum_{\vec i}i_m u(\vec i,t).
\]

Proceed to generating functions $m(\vec x,\vec I,t)=\sum_{\vec i} \prod_k x_k^{i_k} M_{\vec i,\vec I}(t)$.
Note that $\frac{\partial}{\partial x_k}m(\vec x,\vec I,t)=\sum_{\vec i,i_k>0}kx_k^{i_k-1}\prod_{j\neq k}x_j^{i_j} M_{\vec i,\vec I}(t)$ and 
$m(\vec x,\vec I,0)=P(\vec I)$. 

Multiply the master equation (\ref{master_eq3}) by $\prod x_k^{i_k}$ and sum up over all $\vec i$. Thus, obtain the following system
of PDEs
\begin{equation}
	\label{pde_3}
	\begin{array}{l}
\frac{\partial m(\vec x,\vec I,t)}{\partial t}=((\vec x-1)\odot (m(\vec x,\vec I,t)\vec I-\grad_x m(\vec x,\vec I,t)\odot \vec x))^T W(\vec\nu_1-\vec\mu(t)),\\
m(\vec x,\vec I,0)=P(\vec I).
	\end{array}
\end{equation}

Differentiating (\ref{pde_3}) by $x_k$ and summing up over $\vec I$ at all $\vec x=1$, we get ODE for $\vec\mu(t)$
\begin{equation}\label{ode_mu3}
	\begin{array}{l}
\frac{\partial\mu(t)}{\partial t}=(\vec\nu_1-\vec\mu)\odot W(\vec\nu_1-\vec\mu(t)),\\
\vec\mu(0)=0
	\end{array}
\end{equation}

Suppose, that we have a solution (it can be easily obtained numerically) to the system (\ref{ode_mu3}). Assume that
\begin{equation}\label{suggest_m3}
m(\vec x,\vec I,t)=P(\vec I)\prod_k
\left(\frac{1+A_k(t)t}{1+A_k(t)tx_k}\right)^{-I_k}.
\end{equation}

Divide equations in (\ref{pde_3}) by $m(\vec x,\vec I,t)$:
\begin{equation*}
	\begin{array}{l}
\frac{\partial\log m(\vec x,\vec I,t)}{\partial t}=
(\vec x-1)\odot ((\vec I-\grad_x m(\vec x,\vec I,t)\odot\vec x))^T W(\vec\nu_1-\vec\mu(t)).
\end{array}
\end{equation*}
Substitute $\log m(\vec x,\vec I,t)$ from (\ref{suggest_m3}) and equal terms depending on $\vec x$. Namely, we get the following system of ODEs
\begin{equation*}
%	\begin{array}{l}
\vec A'(t)t+\vec A(t)=(\vec 1+\vec A(t)t)\odot W(\vec\nu_1-\vec\mu(t))
%	\end{array}
\end{equation*}
with initial condition $\vec A(0)=\vec\nu_1$.

Thus, 
\begin{equation}\label{abc_3}
%\begin{array}{l}
\vec A(t)=\frac{\exp\{\int_0^t W(\vec\nu_1-\vec\mu(u))du\}-\vec 1}{t},
%\end{array}
\end{equation}
where exponent is assumed to be applied componentwise along with integration.

Inverting generating functions in (\ref{suggest_m3}) and denoting
$p_k(t)=\frac{tA_k(t)}{1+A_k(t)}$,
we obtain the following expression
\begin{equation}
	\label{M_3}
	\begin{array}{l}
M_{\vec i,\vec I}(t)=P(\vec I)\prod_k\binom{I_k}{i_k}p_k(t)^{i_k}(1-p_k(t))^{I_k-i_k}.
	\end{array}
\end{equation}

Consequently,
\begin{equation}
	\label{u_3}
	\begin{array}{l}
		u(\vec i,t)=\sum_{\vec I}P(\vec I)\prod_k\binom{I_k}{i_k}p_k(t)^{i_k}(1-p_k(t))^{I_k-i_k}.
	\end{array}
\end{equation}

%Similarly to (\ref{mu_mn_e}), we have $\vec\mu(t)=\vec p(t)\odot\vec\nu_1$. Now proceed to component sizes calculation.

Define $	\mu_{mn}(t)=\sum_{\vec i}i_mi_n u(\vec i,t)$. 
Quick calculations give relations similar to (\ref{mu_mn_e}):
\begin{equation}\label{mu_mnk_e}
	\begin{array}{l}
	\vec\mu(t)=\vec p(t)\odot\vec\nu_1,\\
	\mu_{mm}(t)=p_m(t)^2\nu_{2,mm}+p_m(t)(1-p_m(t))\nu_m,\\
	\mu_{nm}(t)=p_m(t)p_n(t)\nu_{2,mn}, m\neq n.
%\mu_{000}(t)=1, \mu_{100}(t)=p_A(t)\nu_{100}, \mu_{010}(t)=p_B(t)\nu_{010}, \mu_{001}(t)=p_C(t)\nu_{001},\\
%\mu_{200}(t)=p_A(t)^2\nu_{200}+p_A(t)(1-p_A(t))\nu_{100}, \\
%\mu_{020}(t)=p_B(t)^2\nu_{020}+p_B(t)(1-p_B(t))\nu_{010}, \\
%\mu_{002}(t)=p_C(t)^2\nu_{002}+p_C(t)(1-p_C(t))\nu_{001}, \\
%\mu_{110}(t)=p_A(t)p_B(t)\nu_{110}, \mu_{011}(t)=p_B(t)p_C(t)\nu_{011}, \\
%\mu_{101}(t)=p_A(t)p_C(t)\nu_{010},
%\mu_{111}(t)=p_A(t)p_B(t)p_C(t)\nu_{111}.
	\end{array}
\end{equation}
%
%Unfortunatelly, generated network has three types of edges: in, out, and bidirectional, thus in order to calculate weight-average molecular weight, we have to generalize results, obtained in Kryven \cite{kryven_17} and \cite{newman_01}.
Now proceed to component sizes calculation.

Let 
\[
U(\vec x,t)=\sum_{\vec i}u(\vec i,t)\prod_k x_k^{i_k}
\]
be a GF for the degree distribution. Further we'll omit time dependence for brevity. GFs for the excess distributions (for different types of connections) are defined as 
\begin{equation*}
U_k(\vec x)=\frac{1}{\mu_k}\frac{\partial U(\vec x)}{\partial x_k}
\end{equation*}

Denote $q(s)$ probability that a randomly chosen monomer belongs to a component of size $s$. Corresponding GF is defined as $Q(x)=\sum_s q(s)x^s$. Similarly to (\ref{W_in_bi}), we have the following system of functional equations:
\begin{equation}
	\label{W3}
	\begin{array}{l}
		Q(x)=xU(K\vec Q(x)),
		Q_k(x)=xU_k(K\vec Q(x)),
	\end{array}
\end{equation}
where 
\[
\begin{array}{l}
K_{ij}=\frac{W_{ij}/\mu_j}{\sum_{j}W_{ij}/\mu_j}.
\end{array}
\]
Changing variables $\vec R(x)=K\vec Q(x)$, we get
\begin{equation}\label{W4}
  Q(x)=xU(\vec R(x)), K\vec Q(x)=x K\vec U(\vec R(x)).
\end{equation}
Note that due to normalization $\vec R(1)=\vec 1$.

Differentiating (\ref{W4}) at $x=1$, we get
\begin{equation}
	\label{W3_diff}
	\begin{array}{l}
		Q'(1)=1+\sum_k U_k(\vec R(1))\frac{\partial R_k(1)}{x},\\
		\frac{\partial R_k(1)}{\partial x}=1+\sum_{j,m} K_{km} \frac{\partial U_m(\vec R(1))}{\partial x_j}\frac{\partial R_j(1)}{\partial x}
	\end{array}
\end{equation}

The system of linear equations (\ref{W3_diff}) with respect to variables $\frac{\partial R_j(1)}{\partial x}$ degenerates once determinant of its matrix equals zero. This condition is used to get the gel transition point. It remains to note that derivatives $\frac{\partial U_k(\vec 1)}{\partial x_j}$ can be expressed through $\mu_k$ similarly to (\ref{diff_u_through_m}) using the relations $\mu_k=\left.\left(x_k\frac{\partial}{\partial x_k}\right)U(\vec x)\right|_{\vec x=\vec 1}$ and $\mu_{mn}=\left.\left(x_m\frac{\partial}{\partial x_m}\right)\left(x_n\frac{\partial}{\partial x_n}\right)U(\vec x)\right|_{\vec x=\vec 1}$. Namely, 
\[
\begin{array}{l}
\frac{\partial U_k(\vec 1)}{\partial x_k}=\frac{\mu_{kk}}{\mu_k}-1,\\
\frac{\partial U_m(\vec 1)}{\partial x_n}=\frac{\mu_{mn}}{\mu_m}, m\neq n.\\
\end{array}
\]
Substituting (\ref{mu_mnk_e}), we get
\[
\begin{array}{l}
	\frac{\partial U_k(\vec 1)}{\partial x_k}=p_k(\frac{\nu_{2,kk}}{\nu_k}-1),\\
	\frac{\partial U_m(\vec 1)}{\partial x_n}=\frac{p_{n}\nu_{2,mn}}{\nu_m}, m\neq n.\\
\end{array}
\]

Thus, we obtain the following 
\begin{theorem}
	For systems with several functional groups and matrix of possible connections gel transition point is reached once determinant of the matrix
	\[\tiny
	\left(
	\begin{array}{llll}
		1+K_{11}p_1-\sum\limits_m \frac{K_{1m}p_1\nu_{2,m1}}{\nu_m} & 
		K_{12}p_2-\sum\limits_m \frac{K_{1m}p_2\nu_{2,m2}}{\nu_m}& \cdots & K_{1n}p_n-\sum\limits_m \frac{K_{1m}p_n\nu_{n,mn}}{\nu_m}		\\
		K_{21}p_1-\sum\limits_m \frac{K_{2m}p_1\nu_{2,m1}}{\nu_m} & 
1+K_{22}p_2-\sum\limits_m \frac{K_{2m}p_2\nu_{2,m2}}{\nu_m}& \cdots	& K_{2n}p_n-\sum\limits_m \frac{K_{2m}p_n\nu_{n,mn}}{\nu_m}	\\
\cdots&\cdots&\cdots&\cdots\\
		K_{n1}p_1-\sum\limits_m \frac{K_{nm}p_1\nu_{n,m1}}{\nu_m} & 		K_{n2}p_2-\sum\limits_m \frac{K_{nm}p_2\nu_{n,m2}}{\nu_m} & \cdots &
1+K_{nn}p_n-\sum\limits_m \frac{K_{nm}p_n\nu_{n,mn}}{\nu_m}\\

	\end{array}
\right)
	\]
	equals zero.
\end{theorem}

\section{Appendix}

Here we present numerical results obtained for different monomer configurations.

\subsection{Single monomer type}

\begin{figure}[H]
	\centering
	\includegraphics[width=0.7\linewidth]{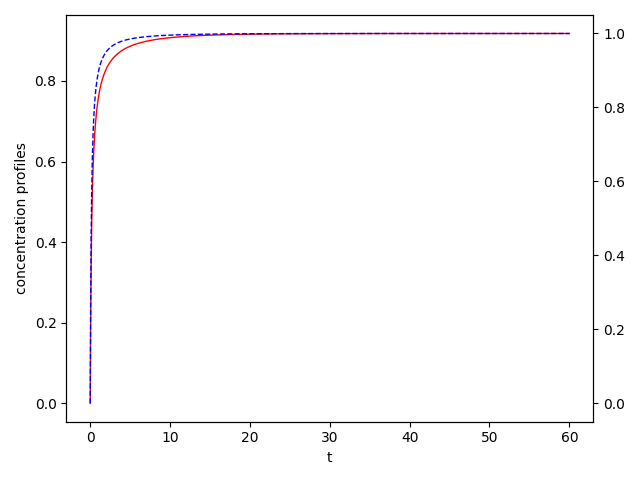}
	\caption[Bond concentration profiles]{Bond concentration profiles}
	\label{fig:singl1}
\end{figure}

Figure \ref{fig:singl1} shows $p_A(t)$ and $p_B(t)$ for monomer configuration $P(I,J)=1, I=2, J=5$. Note that these bond concentration profiles are not proportional (as we expect for purely directed or undirected graph models), however they both converge to $1$. 

Note that if we change reaction rate for B-B type bonds slower $10$ times, the model shows more difference in concentration profiles, but still they both converge to the same value as $t$ goes to infinity.

\begin{figure}[H]
	\centering
	\includegraphics[width=0.7\linewidth]{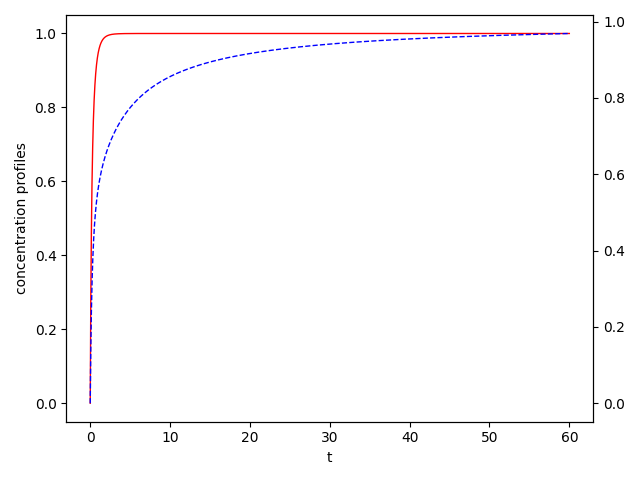}
	\caption[Bond concentration profiles]{Bond concentration profiles with B-B group reacting slower}
	\label{fig:singl2}
\end{figure}

\subsection{Two monomer types configurations}

In this section we present numerical results corresponding to monomers with functional groups of two kinds. Let us start with fully oriented network case when we have homofunctional polycondensation of monomers with 2 and 5 functional groups correspondingly. This case can be easily derived from the paper \cite{Kryven2}. Note that we have $\mu_1(t)\equiv\mu_2(t)$. Since step-growth polymers increase in molecular weight at a very slow rate, we're interested both in conversion rate at gel transition point and (relative) time.

\begin{figure}[h]
	\centering
	\includegraphics[width=0.7\linewidth]{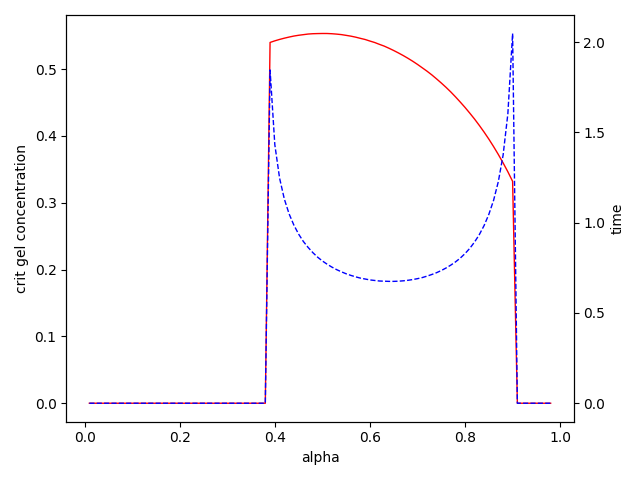}
	\caption{Gel transition point for different stoichiometric ratios}
	\label{fig:p2-5-noself}
\end{figure}

From the modelling results in Fig.~\ref{fig:p2-5-noself} it can be easily seen that largest conversion rate and smallest time are attained in different points, however from practical point of view stoichiometric ratios (alpha) in quite a wide range produce satisfactory results.

Also, we can point out that we have two intervals for ratios where transition is not observed in reasonable time. Obsiously, this is due to the lack of monomers of one or another type.

\subsection{Self-polymerizing monomers}

Now we assume that the other is a self-polymerizing monomer with 5 functional groups. Note that corresponding network contains both directed (which indicate A-B bonds) and undirected edges (indicating B-B bonds). 

\begin{figure}[h]
	\centering
	\includegraphics[width=0.7\linewidth]{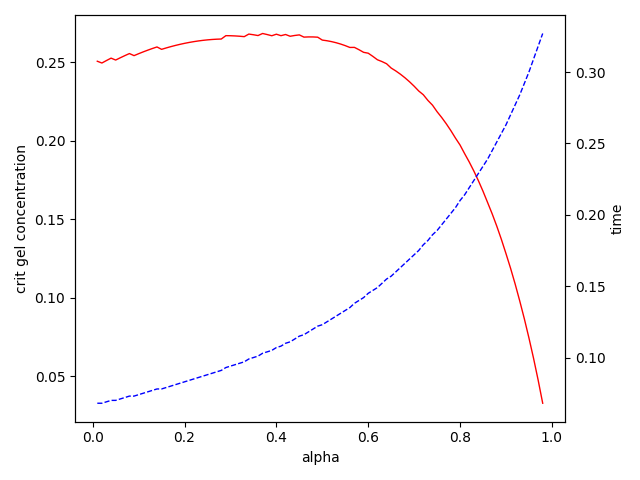}
	\caption{Gel transition point for different stoichiometric ratios}
	\label{fig:p2-5-w1}
\end{figure}

Fig.~\ref{fig:p2-5-w1} shows quite natural behavior of the conversion rate depending on the ratio between monomers. Obviously, if we have lots of A-A monomers, gel fraction is going to be low and it is going to take much more time before the large component appears. Trying to find optimal (e.g. showing largest conversion rate and lowest time) ratio gives us slightly smaller values compared to the heterofunctional polycondensation.

\subsection{Accounting for copolymerization reaction rate}

Now let us see how copolymerization reaction constants affect gel transition for the above example. 

\begin{figure}[H]
	\centering
	\includegraphics[width=0.7\linewidth]{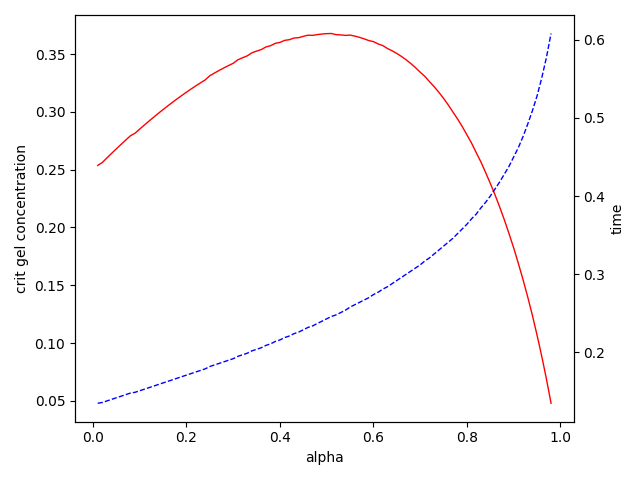}
	\caption{Gel transition point with twice lower rate of B-B reaction}
	\label{fig:p2-5-w5}
\end{figure}

In the Fig.~\ref{fig:p2-5-w5} we have results for the case when self-polymerization rate is twice less probable than A-B reaction. Note that compared to the original model, optimal range of ratios shifts towards the one in heterofunctional polycondensation case. 

If we take self-polymerization rate even lower, we observe that plot converges to the one in Fig.~\ref{fig:p2-5-noself} which indicates self-consistency of the model.

\begin{figure}[h]
	\centering
	\includegraphics[width=0.7\linewidth]{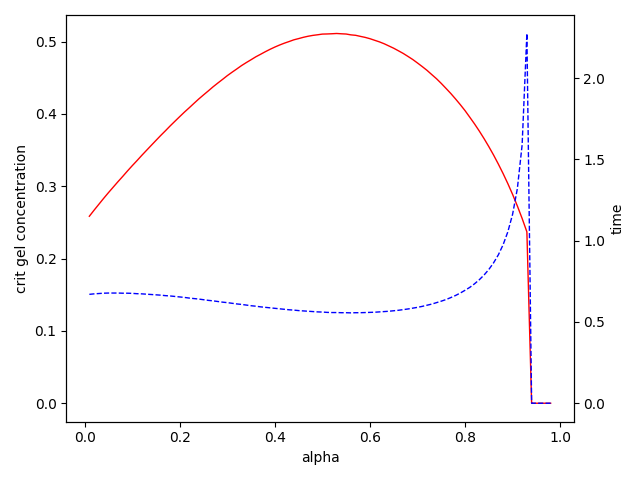}
	\caption{Gel transition point with 10 times less probable B-B bond}
	\label{fig:p2-5-w_1}
\end{figure}

\subsection{Systems with multiple functional group types}

In this section we show the ability of obtained equations to predict gel transition in even more complicated systems. Here we consider a system with 3 kinds of bonds, heterofunctional A-A, B-B, and monomers with 5 functional groups which are able to connect to A but also have self-polymerization capabilities (though at 10 times lower rate). Now we have two ratios, $\alpha$ and $\beta$ which denote corresponding fractions of monomers of first two types, $\alpha+\beta<1$, and 3rd type monomers are observed at fraction $1-\alpha-\beta$. Thus, we can calculate overall fraction of converted monomers (pCrit), time to the transition point.
Results are presented in the following color plots. Combined plot adds constant time levels to the pCrit plot.

\begin{figure}[H]
	\centering
	\includegraphics[width=0.7\linewidth]{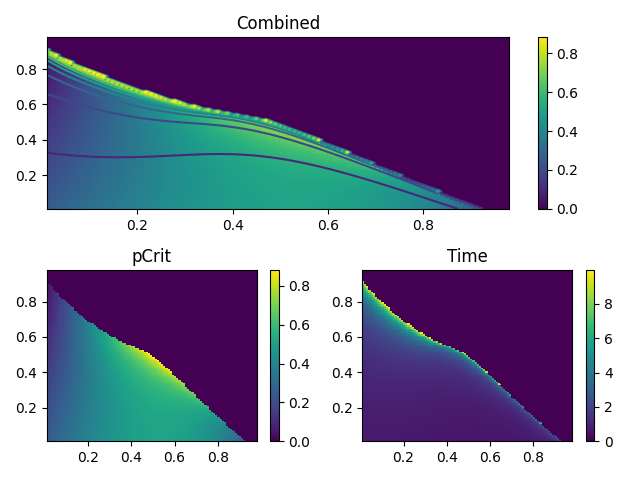}
	\caption{Gel transition point in a complex system}
	\label{fig:p2-2-5}
\end{figure}
 
 Analysis of such figures becomes more complicated, but we can see that in this case it is beneficial to have most of concentration distributed equally between first two kinds of monomers while leaving small fraction for the 3rd component.

\end{document}